\documentclass[12pt]{article}
\usepackage{amsmath}
\usepackage{amsfonts}
\usepackage{latexsym}
\newtheorem{Lem}{Lemma}
\newtheorem{Cor}{Corollary}
\newtheorem{Pro}{Proposition}
\newtheorem{Rem}{Remark}
\newtheorem{The}{Theorem}
\begin{document}
\title{Topology of the Set of Smooth Solutions\\
to the Liouville Equation\\}
\author{Kazimierz Br\c{a}giel$^1$, W{\l}odzimierz Piechocki$^2$}
\maketitle
\hspace{3mm}
\begin{tabbing}
\small $^1$Institute of Physics, Warsaw University Branch in 
Bia{\l}ystok,\\ 
\small Lipowa 41, 15-424 Bia{\l}ystok, Poland. 
E-mail: kbragiel@fuw.edu.pl \\  \small 
$^2$Field Theory Group, So{\l}tan Institute for Nuclear Studies,\\  
\small Ho\.{z}a 69, 00-681 Warsaw, Poland. E-mail: piech@fuw.edu.pl
\end{tabbing}
\begin{abstract}
We prove that the space of smooth initial data and the set of smooth
solutions of the Liouville equation are homeomorphic.
\end{abstract}
\pagebreak

\section{Introduction}

There is a considerable need for the method which can be used
to the quantization of \emph{nonlinear} field theories. We still do 
not know how to quantize, in mathematically correct way, the Einstein 
theory of gravitation and the Yang-Mills theories. In the case of 
simple systems with \emph{finitely} many degrees of freedom the method 
of geometric quantization ( see, e.g.,~\cite{5} ) seems to be 
satisfactory. It is not clear, however, if we can generalize this 
method to the case of nonlinear field theories. The main problem is 
that such theories have \emph{infinitely} many degrees of freedom and 
the set of solutions to the field equations \emph{cannot be} 
a vector space.
The 2-dim Liouville field equation~\cite{2}
\[ \left( \partial^2_t - \partial^2_x \right) F(t,x) + \frac{m^2}{2}
\exp F(t,x) = 0,\;\;\;\;\;m>0  \]
is a simple model that we are going to use for verification of ideas 
connected 
with generalization of the geometric quantization method to nonlinear 
field theories.

In the geometric quantization procedure one assumes that 
the phase space of a given classical theory is a \emph{manifold}. 
In the case of simple systems with n-degrees of freedom 
one can easily identify the phase space to be a manifold modelled on
$\mathbf{R}^{2n}$. In the case of field theory we would like to have 
an object which we could call a manifold modelled on a Fr{\'{e}}chet 
space. The aim of this paper is to show that the space of smooth 
initial data and the set of smooth solutions to the Liouville 
equation are \emph{homeomorphic}. Therefore, the set of solutions 
has the structure of a topological manifold modelled on the space of
initial data.

It follows from the existence of homeomorphism that 
if there exists a series of 
smooth initial data convergent almost uniformly with all its 
derivatives to $(f,g) \in C^{\infty}(\mathbf{R}) \times C^{\infty}
(\mathbf{R})$, then the series of corresponding solutions of the 
Liouville equation converges almost uniformly with all its derivatives 
to the solution corresponding to the initial data $(f,g)$.\\
$\mathbf{Notation:}$ Through the paper we use the spaces
$C^{\infty}(\mathbf{R}^M, \mathbf{R}^N)$,\\ where $ M, N \in \mathbf{N}
^{\times};\;\; 
\mathbf{N}^{\times} := \mathbf{N}\setminus \{0\},\;\;\mathbf{N} =  
\{0,1 \ldots \}$;\\ for $\alpha \in \mathbf{N}^{\times},
\;K_{\alpha} := [\;- \alpha, 
\alpha\;],$ thus $\mathbf{R}^M = \bigcup_{\alpha = 1}^{\infty} K^M_
{\alpha}$.\\ If $\beta = (\beta_1, \ldots ,\beta_M) \in 
\mathbf{N}^M$, 
then $\partial^{\beta} := \partial^{\beta_1}_1 \cdot \partial^{\beta_2}
_2 \cdot \ldots \cdot \partial^{\beta_M}_M,\;\;\mid \beta \mid := 
\sum_{j=1}^{M} \beta_{j}$.
\newpage
\noindent
Topology on $C^{\infty}(\mathbf{R}^M, 
\mathbf{R}^N)$ is given by the family 
of seminorms \\$\left( p^{MN}_
{\alpha \beta}\right)_{{(\alpha, \beta)} \in \mathbf{N}^{\times}\times 
\mathbf{N}^{M}}$ defined by
\begin{displaymath}
p^{MN}_{\alpha \beta} : C^{\infty}(\mathbf{R}^M, \mathbf{R}^N) \ni 
f \longrightarrow p^{M N}_{\alpha \beta}(f) := 
 \underset{x \in K_{\alpha}^M}{\sup}
\parallel (\partial^{\beta} f)(x)\parallel \in \mathbf{R},
\end{displaymath}
where
\[\parallel \cdot \parallel \;: \mathbf{R}^N \ni y \rightarrow 
\parallel y \parallel := \left(\sum_{j=1}^{N}(y^j)^2 \right)^{\frac{1}
{2}} \in \mathbf{R}. \]
Chapter III of ~\cite{4} is our source of information 
on topological vector 
spaces ( in particular Fr{\'{e}}chet spaces ).\\ 
To denote  maps 
$\mathbf{R} \rightarrow \mathbf{R}$ we use small Latin letters,\\ 
$(f,g,h, \ldots \in C^{\infty}(\mathbf{R}) := C^{\infty}(\mathbf{R},
\mathbf{R}))$; for maps $\mathbf{R}^2 \rightarrow \mathbf{R}$ capital 
Latin letters are used, $(F,G,H, \ldots \in C^{\infty}(\mathbf{R}^2))$;
 $\mathbf{R} \rightarrow \mathbf{R}^2$ maps are denoted by  capital 
 Greek letters, $(\Phi, \Psi, \Omega, \ldots \in C^{\infty}(\mathbf{R},
 \mathbf{R}^2))$. To shorten notation, we write $p_{\alpha \beta} := 
p^{11}_{\alpha \beta},\;\; r_{\alpha \beta} := p^{12}_{\alpha \beta},\;
\; q_{\alpha \beta} := p^{21}_{\alpha \beta}$ for $\; \beta = (
\beta_1, \beta_2); \;\; \partial_t =\partial_1,\;\;\partial_x = 
\partial_2\;$ for $\;(t,x) \in \mathbf{R}^2$. A norm of a linear map 
$A \in B(\mathbf{R}^2) := B(\mathbf{R}^2, \mathbf{R}^2)$ is defined by 
\[ \parallel \cdot \parallel \;: B(\mathbf{R}^2) \ni A \rightarrow 
\parallel A \parallel := 
\underset{\parallel x \parallel =1}{\sup}
\parallel A(x) \parallel \in \mathbf{R}.  \]
$ C^{\infty}_{+}(\mathbf{R}^2) := \{F \in C^{\infty}(\mathbf{R}^2) 
: F(\mathbf{R}^2) 
\subseteq \;]\;0, \infty\;[\;\}$ is a topological space with topology 
induced 
from $C^{\infty}(\mathbf{R}^2)$. From now on, $\;\Box := \partial^2_t-
\partial^2_x,\\ \mathcal{M} := \{F \in C^{\infty}(\mathbf{R}^2) : 
\Box F = -(m^2/2) \exp F\},\;\;$ where $ m > 0$ 
is a fixed real number.

In Sec.(2) we quote some results of ~\cite{1} concerning the solution 
of the Cauchy problem for the Liouville equation and we prove these 
results. The proof that the space of smooth initial data and the set 
of smooth solutions are homeomorphic is given in Sec.(3). We make some 
remarks in the last Section.

\section{Smooth Solutions of the Liouville Equation}

We examine properties of a smooth solution to the Liouville 
equation, i.e., of class $C^{\infty}(\mathcal{O})$, where $\mathbf{R}^2
\supseteq \mathcal{O}$ is an open subset different to $\emptyset$ . 
However, all 
proofs can be easily modified to include the solutions of class 
$C^k(\mathbf{R}^2),\;$ for $k\geq 2.$
\begin{Lem}
Let $g_i \in C^{\infty}(\mathbf{R})\;$ for $i \in \overline{1,4}$ are 
such that
\begin{equation}
g_1 g_3^{\prime} - g_1^{\prime} g_3 = 1,\;\;\;g_2 g_4^{\prime} - 
g_2^{\prime} g_4 = - 1
\end{equation}
and let 
\[ G : \mathbf{R}^2 \ni (t,x) \rightarrow G(t,x) := g_1(x + t) g_2
(x-t) + g_3 (x+t)g_4(x-t) \in \mathbf{R}. \]
Then we have
\[ \left(G^{-1}(0) \neq \emptyset \right) \Longrightarrow \left( 
G^{-1}(0)
\cap \{(0,x) \in \mathbf{R}^2 : x \in \mathbf{R}\} \neq \emptyset
\right). \]
\end{Lem}
\emph{Proof}. Let $\xi := x+t,\;\; \eta := x-t$. Making use of
\begin{equation}
\left( G(t_0,x_0) = 0 \right) \Longleftrightarrow \left( g_1(\xi_0) g_2
(\eta_0) = - g_3(\xi_0) g_4(\eta_0) \right)
\end{equation}
we see that the condition
\[ \left\{ \begin{array}{l}G(t_0,x_0) = 0 \\ (\partial _t G)(t_0,x_0) 
= (\partial _x G)(t_0,x_0)
\end{array}
\right. \]
leads to
\[ g_1(\xi_0) g_2^{\prime}(\eta_0) + g_3(\xi_0) g_4^{\prime}(\eta_0) 
= 0. \]
Multiplying this equation by $\;g_2(\eta_0)\;$ and using (1) gives 
$\;g_3(\xi_0) = 0$. Similarly, multiplying by $\;g_4(\eta_0)\;$ 
leads to $\;g_1(\xi_0)=0$.
Thus, we have
\[ \left( \begin{array}{lcl} G(t_0,x_0)& =& 0 \\ (\partial_t G)
(t_0,x_0)& 
=& (\partial_x G)(t_0,x_0) \end{array} \right) \Longrightarrow
\left( g_1(t_0,x_0) = 0 = g_3(t_0,x_0) \right), \]
contrary to (1).
In the same manner we can see that
\[ \left( \begin{array}{lcl} G(t_0,x_0) & =& 0 \\ (\partial_tG)
(t_0,x_0) &=& - (\partial_xG)(t_0,x_0) \end{array} \right) 
\Longrightarrow
(g_2(t_0,x_0) = 0 = g_4(t_0,x_0) ),  \]
which again contradicts (1).
Therefore, we have
\begin{equation}
\left( (t_0,x_0) \in G^{-1}(0) \right) \Longrightarrow \left( 
(\partial_t G)(t_0,x_0) \neq \pm (\partial_x G)(t_0,x_0) \right).
\end{equation}
The condition $\;(\partial_x G)(t_0,x_0) = 0\;$ means that
\begin{equation}
0 = g_1^{\prime}(\xi_0)g_2(\eta_0) + g_1(\xi_0)g_2^{\prime}(\eta_0) + 
g_3^{\prime}(\xi_0)g_4(\eta_0) + g_3(\xi_0)g_4^{\prime}(\eta_0).
\end{equation}
Multiplying (4) by $\;g_1(\xi_0)g_4(\eta_0)\;$ and using (2) gives
\begin{equation}
g_1(\xi_0)^{2} + g_4(\eta_0)^2 = 0.
\end{equation}
Similarly, multiplying (4) by $\;g_3(\xi_0)g_2(\eta_0)\;$ and using (2) 
leads to 
\begin{equation}
g_3(\xi_0)^2 + g_2(\eta_0)^2 = 0.
\end{equation}
Since (5) and (6) contradict (1), we conclude that
\begin{equation}
\left( (t_0,x_0) \in G^{-1}(0) \right) \Longrightarrow \left(
(\partial_xG)(t_0,x_0) \neq 0 \right),
\end{equation}
which means that either $\;G^{-1}(0) = \emptyset\;$ or $\;G^{-1}(0)\;$  
is a one-dimensional $C^{\infty}$ submanifold of $\mathbf{R}^2$. 
Suppose that $\;G^{-1}(0) \neq \emptyset\;$ and let us denote by $M$ 
the 
connected component of $\;G^{-1}(0)$. By (3) we have that $M$ cannot be 
a compact subset of $\mathbf{R}^2$. Since $M$ is closed in $\mathbf
{R}^2$, it cannot be bounded in $\mathbf{R}^2$. From (7) we conclude 
that $M$ can be parametrized by $t \in \mathbf{R}$. Let $\pi
: \mathbf{R}^2 \ni (t,x) \rightarrow \pi (t,x) := t \in 
\mathbf{R}$. Since $M$ is closed, $\pi (M)\subseteq \mathbf
{R}\;$ is closed in $\mathbf{R}$. Hence $\emptyset \neq \pi (M)
\subseteq \mathbf{R}\;$ is both closed and open (homeomorphic to 
$\mathbf{R}$). Therefore $\pi (M) = \mathbf{R}$. Finally, 
we obtain
\[ G^{-1}(0) \cap \{(0,x) \in \mathbf{R}^2 : x \in \mathbf{R} \} 
\neq \emptyset. \;\;\;\;\; \Box \]
By ~\cite{1} we get
\begin{Lem}
Let $\mathbf{R}^2\supseteq \mathcal{O}$ be an open subset and let
\[ \mathcal{O}_1 := \{(x +t) \in \mathbf{R} \;:\; (t,x) \in \mathcal{O}
\},\;\;\; \mathcal{O}_2 := \{(x-t) \in \mathbf{R}\;:\;(t,x) \in 
\mathcal{O}  \}. \]
Suppose $ F \in C^{\infty}(\mathcal{O})$, then the following are 
equivalent:
\begin{enumerate}
\item $\;\;\Box F = -\frac{m^2}{2} \exp F.$
\item \,\,\,There exist $\;g_1,g_3 \in C^{\infty}(\mathcal{O}_1)\;$ and 
$\;g_2,g_4 \in C^{\infty}(\mathcal{O}_2)\;$ satisfying
$ \;g_1 g_3^{\prime} - g_1^{\prime} g_3 = 1 \;\;$and
 $\;\;  g_2 g_4^{\prime} - g_2^{\prime} g_4 = - 1 $
such that
$\; F : \mathcal{O}\ni (t,x) \rightarrow F(t,x) := -\log 
\frac{m^2}{16}\left[ 
g_1(x+t)g_2(x-t) + g_3(x+t)g_4(x-t) \right]^2 \in \mathbf{R}.   $
\end{enumerate}
\end{Lem}
Lemmas (1) and (2) lead to
\begin{Cor}
If $\;\mathbf{R}^2 \supseteq \mathcal{O}\; $ is an open set, such that 
$\;\{(0,x) \in \mathbf{R}^2\;:\;x \in \mathbf{R}\} \subseteq \mathcal{O}
\;$ and 
$\;F \in C^{\infty}(\mathcal{O})\;$ satisfies the Liouville equation 
$\;\Box F = - (m^2/2) \exp F,$ then there exists $\;\;\widetilde{F}
\in C^{\infty}
(\mathbf{R}^2)\;$ such that $\;\Box \widetilde{F} = -( m^2/2) \exp 
\widetilde{F}$ 
and $\widetilde{F}_{\mid _{\mathcal{O}}} = F.$
\end{Cor}
\begin{Pro}
Suppose $\;f_1,f_2,g_1,g_2,g_3,g_4 \in C^{\infty}(\mathbf{R})\;$ and let
\begin{equation}
g_1g_3^{\prime}-g_1^{\prime}g_3=1,
\end{equation}
\begin{equation}
g_2g_4^{\prime}-g_2^{\prime}g_4=-1,
\end{equation}
\begin{eqnarray}
G:\mathbf{R}^2 \ni (t,x)\rightarrow G(t,x)& := & g_1(x+t)g_2(x-t)
+ \nonumber \\ & & 
g_3(x+t)g_4(x-t) \in \mathbf{R},
\end{eqnarray}
\begin{equation}
u := \frac{1}{16} \left[ (f_1^{\prime}-f_2)^2-4 (f_1^{\prime}-f_2)
^{\prime} + m^2 \exp f_1 \right],
\end{equation}
\begin{equation}
w := \frac{1}{16}\left[ (f_1^{\prime}+f_2)^2 - 4(f_1^{\prime}+f_2)
^{\prime} + m^2 \exp f_1 \right].
\end{equation}
Then
\begin{enumerate}
\item If $\;g_1,\ldots ,g_4\;$ satisfy the equations
\begin{equation}
g_i^{\prime \prime}=u g_i\;\;\;\;\;\mbox{for}\;\;\;i=2,4 
\end{equation}
and
\begin{equation}
g_j^{\prime \prime}=w g_j\;\;\;\;\;\mbox{for}\;\;\;j=1,3
\end{equation}
then the map
\begin{equation}
F:\mathbf{R}^2 \ni (t,x) \rightarrow F(t,x) := - \log \frac {m^2}{16}
G^2(t,x) \in \mathbf{R}
\end{equation}
is a solution of the Liouville equation with the initial data 
\begin{equation}
\left\{ \begin{array}{ccl}F(0,\cdot)& =& f_1 \\
(\partial_t F)(0,\cdot)& = & f_2.
\end{array} \right. 
\end{equation}
\item A solution of the Liouville equation satisfying (16) is given 
by (15), where $\;g_1,\ldots,g_4\;$ satisfy (13) and (14).
\end{enumerate}
\end{Pro}
\emph{Proof (1.)} Suppose $\;g_1,g_3\;$ and $\;\widetilde{g}_1, 
\widetilde{g}_3\;$ 
satisfy (8) and (14). Hence there exists $\left( \begin{array}{cc}
a&b\\c&d \end{array} \right) \in SL(2,\mathbf{R})\;$ such that
$\;g_1 = a\widetilde{g}_1 + b \widetilde{g}_3$ and $ g_3 =c \widetilde{g
}_1 + d \widetilde{g}_3.\;$ Such an exchange of functions corresponds  
to the Bianchi transformation~\cite{1} and does not change the form of 
solution (15). If $g_2$ and $g_4$ satisfy (9) and (13), 
then the functions~\cite{1} 
\begin{equation}
g_1 := - \frac {4}{m}\exp \left( -\frac {1}{2} f_1(\cdot)\right)
\left[ g_4^{\prime} + \frac{1}{4}(f_1^{\prime}-f_2)g_4 \right ], 
\end{equation}
\begin{equation}
g_3 := \frac{4}{m} \exp \left(- \frac{1}{2}f_1(\cdot) \right) 
\left [g_2^{\prime} + \frac{1}{4}(f_1^{\prime}-f_2)g_2 \right]
\end{equation}
satisfy (8) and (14). In what follows we assume that the general 
form of $g_1$ and $g_3$ are given by (17) and (18). One can easily 
check that for $G$ defined by (10) we have
\[ \forall x \in \mathbf{R}\;:\;G(0,x) = \frac{4}{m}\exp \left(
- \frac {1}{2}f_1(\cdot) \right) > 0. \]
By Lemma 1 we have that $F$ given by (15) is well defined on $\mathbf
{R}^2 \\(F \in C^{\infty}(\mathbf{R}^2))$ and satisfies (16).\\
\emph{(2.)} 
We have shown that there exists $F \in C^{\infty}(\mathbf{R}
^2)$ satisfying (16). By Lemma 2, $F$ is of the form (15) for $\;g_1,
\ldots,g_4\in C^{\infty} (\mathbf{R})\;$ satisfying (8) and (9). Now, 
we shall show (see~\cite{1}) that $g_1, \ldots, g_4$ can be a solution 
of (13) and (14) with $u$ and $w$ given by (11) and (12). \\
Let
\[ \aleph : \mathbf{R} \ni x \rightarrow \aleph (x) := G(0,x) \in 
\mathbf{R} \]
\[ \hbar : \mathbf{R} \ni x \rightarrow \hbar (x) := g_1^{\prime}(x)
g_2^{\prime}(x) + g_3^{\prime} (x)g_4^{\prime}(x) \in \mathbf{R}. \]
Eq.(16) means that
\[ \forall x \in \mathbf{R} \;\;\;\aleph (x) \neq 0 \]
and
\begin{equation}
\aleph = \epsilon \frac{4}{m} \exp \left(-\frac{1}{2}f_1(\cdot) 
\right),
\end{equation}
where
\[ \epsilon := \left\{\begin{array}{rl}1 &\mbox{for}\;\; \aleph > 0\\
-1 & \mbox{for}\;\; \aleph < 0 \end{array} \right..\]
By (19) we get $\;\aleph\;^{\prime}/\aleph = -\frac{1}{2} f_1^{\prime}
\;$ and
\begin{eqnarray}
f_2 &=&(\partial_t F)(0,\cdot) = \frac{2}{\aleph}\left(\aleph\;^{\prime}
-2 (g_1^{\prime}g_2 + g_3^{\prime}g_4) \right) = - f_1^{\prime}-
\frac{4}{\aleph}(g_1^{\prime}g_2 + g_3^{\prime}g_4) \nonumber\\
 &=& \frac{2}{\aleph}\left( 2(g_1g_2^{\prime}+ g_3^{\prime}g_4)-
 \aleph\;^{\prime} \right)
 = f_1^{\prime}+\frac{4}{\aleph}(g_1g_2^{\prime}+g_3g_4^{\prime})
 \nonumber
 \end{eqnarray}
which leads to
\begin{equation}
\frac{1}{4}(f_1^{\prime}+f_2) \aleph = -g_1^{\prime}g_2-g_3^{\prime}
g_4,
\end{equation}
\begin{equation}
\frac{1}{4}(f_1^{\prime}-f_2)\aleph = g_1g_2^{\prime}+g_3g_4
^{\prime}.
\end{equation}
Suppose now that $\;g_1,g_3\;$ satisfy (13) and $\;g_2,g_4\;$ 
satisfy (14) 
for some $\;u,w \in C^{\infty}(\mathbf{R})$. Taking derivative of (20) 
yields
\begin{equation}
\frac{1}{4}(f_2+f_1^{\prime})^{\prime}\;\aleph +\frac{1}{4}(f_2+f_1
^{\prime})\;\aleph\;^{\prime}=- w \aleph - \hbar
\end{equation}
By analogy, we get from (21):
\begin{equation}
\frac{1}{4}(f_2 - f_1^{\prime})^{\prime}\;\aleph + \frac{1}{4}(f_2-f_1
^{\prime})\;\aleph\;^{\prime} = u\;\aleph + \hbar
\end{equation}
By direct calculations, we get
\begin{equation}
(\partial_t^2 F)(0,\cdot) = \frac{1}{2} f_2^2 + 4 \frac{\hbar}{\aleph}
- 2(u+w)
\end{equation}
Since $\;(\partial_x^2 F)(0,\cdot)= f_1^{\prime \prime}\;$, we get
\begin{equation}
f_1^{\prime \prime}= 2\left( \frac{\aleph\;^{\prime}}{\aleph}\right)^2
- 2 \frac{\aleph\;^{\prime \prime}}{\aleph}
\end{equation}
The map $F$ satisfies the Liouville equation $\;\Box F=-(m^{2}/2)\exp F
\;$on $\mathbf{R}^2$. By (24) and (25), for $\;(0,x) \in \mathbf{R}^2\; 
$ we get
\begin{equation}
\frac{1}{2}f_2^2+ 4 \frac{\hbar}{\aleph}-2(u+w)+ 2 \frac{\aleph^{\;
\prime \prime}}{\aleph}-2 \left(\frac{\aleph^{\;\prime}}{\aleph}
\right)^2
= - \frac{m^2}{2}\exp f_1.
\end{equation}
Since
\[ 2\;\aleph^{\;\prime \prime}/\aleph = 2(u+w)+ 4\;\hbar/\aleph \;\;\;
\mbox{and}\;\;\left(\aleph^{\;\prime}/\aleph \right)^2 = (1/4)(f_1^
{\prime})^2  \]
Eq. (26) leads to
\begin{equation}
\frac{\hbar}{\aleph}=-\frac{1}{16}\left[ (f_2-f_1^{\prime})(f_2+f_1
^{\prime})-m^2\exp f_1 \right].
\end{equation}
By (27) and (22) we get
\begin{equation}
w=\frac{1}{16}(f_1^{\prime}-f_2)^2 - \frac{1}{4}(f_1^{\prime}-f_2)^
{\prime}+ \frac{m^2}{16}\exp f_1.
\end{equation}
Similarly, (27) and (23) give
\begin{equation}
u= \frac{1}{16}(f_1^{\prime}+f_2)^2 - \frac{1}{4}(f_1^{\prime}+f_2)
^\prime + \frac{m^2}{16} \exp f_1. \;\;\;\;\;    \Box
\end{equation}

\section{Homeomorphism of the Space \\of Initial Data 
	 and the Set of Solutions}

We define the following mappings:
\begin{equation}
\mathcal{A}:C^{\infty}(\mathbf{R})^2 \ni (f_1,f_2) \longrightarrow 
\mathcal{A}(f_1,f_2):= u \in C^{\infty}(\mathbf{R}),
\end{equation}
where $u$ is given by (11).
\begin{equation}
\mathcal{B}:C^{\infty}(\mathbf{R}) \ni u \longrightarrow
\mathcal{B}(u):= (g_2,g_4) \in C^{\infty}(\mathbf{R})^2,
\end{equation}
where $g_2$ and $g_4$ are defined by
\begin{displaymath}
\left\{ \begin{array}{lcl}g_2^{\prime \prime} &=&ug_2\\ g_2(0)&=&0\\
g_2^{\prime}(0)&=&1 \end{array}\right.  \;\;\;\mbox{and}\;\;\;\;
\left\{ \begin{array}{lcl}g_4^{\prime \prime}&=&ug_4\\g_4(0)&=&1\\
g_4^{\prime}(0)&=&0 \end{array} \right.
\end{displaymath}
\begin{Rem}
Maps $g_2$ and $g_4$ satisfy (9).
\end{Rem}
\begin{equation}
\mathcal{C}:C^{\infty}(\mathbf{R})^4 \ni (f_1,f_2,g_2,g_4)
\longrightarrow \mathcal{C}(f_1,f_2,g_2,g_4):=(g_1,g_3) \in
 C^{\infty}(\mathbf{R})^2,
\end{equation}
where $g_1$ and $g_3$ are given by (17) and (18).
\begin{equation}
\mathcal{D}:C^{\infty}(\mathbf{R})^2 \ni (f_1,f_2)\longrightarrow 
\mathcal{D}(f_1,f_2):=(f_1,f_2,f_1,f_2) \in  C^{\infty}(\mathbf{R})^4,
\end{equation}
\begin{equation}
\mathcal{E}:C^{\infty}(\mathbf{R}^2) \ni G \longrightarrow
\mathcal{E}(G):=\frac{m^2}{16}G^2 \in C^{\infty}(\mathbf{R}^2),
\end{equation}
\begin{equation}
\mathcal{N}:C^{\infty}_{+}(\mathbf{R}^2) \ni H \longrightarrow 
\mathcal{N}(H):= -\log H \in C^{\infty}(\mathbf{R}^2),
\end{equation}
\begin{equation}
\mathcal{G}:C^{\infty}(\mathbf{R})^4 \ni (g_1,g_3,g_2,g_4)
\longrightarrow \mathcal{G}(g_1,g_3,g_2,g_4)
:=G \in C^{\infty}(\mathbf{R}^2),
\end{equation}
where $G$ is defined by (10).
\begin{equation}
\mathcal{H}:=\mathcal{E}\circ \mathcal{G}\circ (\mathcal{C}\times 
id_2)\circ (id_2 \times \mathcal{D})\circ (id_2 \times \mathcal{B})
\circ (id_2 \times \mathcal{A})\circ \mathcal{D},
\end{equation}
where $\;id_2:=id_{C^{\infty}(\mathbf{R})^2}$.
\begin{Rem}
\(\;\; \mathcal{H}:C^{\infty}(\mathbf{R})^2 \ni (f_1,f_2)
\longrightarrow 
\mathcal{H}(f_1,f_2):= \frac{m^2}{2}G^2 \in C^{\infty}(\mathbf{R}^2)\).
\end{Rem}
By Lemma 1 we have
\begin{Cor}
\(\;\;\;\; \mathcal{H}\left(C^{\infty}(\mathbf{R})^2\right)\subseteq
C^{\infty}_{+}(\mathbf{R}^2) \)
\end{Cor}
Let
\begin{equation}
\mathcal{I}:=\mathcal{N}\circ \mathcal{H}.
\end{equation}
By Proposition 1 we get
\begin{Cor}
$\;\;\;\mathcal{I}:C^{\infty}(\mathbf{R})^2 \longrightarrow \mathcal
{M}\;\;\;$ is a bijection.
\end{Cor}

Now comes the main theorem.
\begin{The}
The mapping $\;\;\mathcal{I}:C^{\infty}(\mathbf{R})^2 \longrightarrow
\mathcal{M}\;\;$ defined by (38) is a homeomorphism.
\end{The}
Before we give the proof, let us prove a few Lemmas.\\
We define some auxiliary maps and sets.\\
For $\beta \in \mathbf{N}^{\times}$:
\[ \mathcal{R}(\beta):= \{ a \in \mathbf{N}^{\beta}:\sum_{j=1}
^{\beta}j a_j=\beta \}, \]
\[ \mathcal{P}_{\beta}:\mathcal{R}(\beta) \ni a\longrightarrow 
\mathcal{P}_{\beta}(a):=\frac{\beta !}{\prod_{j=1}^{\beta}(j !)
^{a_j}a_j !} \in\mathbf{N}, \]
\[ l_{\beta}:\mathcal{R}(\beta) \ni a \longrightarrow l_{\beta}(a):=
\sum _{j=1}^{\beta} a_j \in \mathbf{N}. \]
For $(\lambda,\mu) \in\mathbf{N}^\times \times\mathbf{N}$:
\[ \mathcal{R}(\lambda,\mu):=\{ a \in \mathbf{N}^{\lambda+1}:
\sum _{j=1}^{\lambda}ja_j=\lambda\;,\;\sum_{j=0}^{\lambda}a_j=\mu\},\]
\[ \mathcal{W}_{\lambda ,\mu}:\mathcal{R}(\lambda ,\mu)\ni a
\longrightarrow \mathcal{W}_{\lambda,\mu}(a):=\frac{\mu !}{a_0}\;
\frac{\lambda !}{\prod _{j=1}^{\lambda}(j!)^{a_j}a_j !} 
\in \mathbf{N}, \]
\[ c:=(c^1,\ldots, c^{\beta})\in \overset{\beta}{\underset{i=1}{X}}
\mathcal{R}(\lambda
_i,\mu_i), \]
\[ \mathcal{T}(\lambda,\mu):=\{a\in\mathbf{N}^{\lambda}:\sum_{j=1}
^{\lambda} a_j=\mu \}, \]
\[ \mathcal{N}_{\lambda,\mu}:\mathcal{T}(\lambda,\mu) \ni a 
\longrightarrow \mathcal{N}_{\lambda,\mu}(a):=\frac{\mu !}{\prod_{j=1}
^{\lambda}a_j !} \in \mathbf{N}.  \]
\begin{Lem}
Let $\;\mathbf{R}\supseteq \mathcal{O}_1\;\;$ and $\;\;\mathbf{R}^2 
\supseteq \mathcal{O}_2 \;$ be some open sets,
\[ h \in C^{\infty}(\mathcal{O}_1),\;\;J \in C^{\infty}
(\mathcal
{O}_2)\;\;\;\mbox{and}\;\;\;\forall(t,x) \in \mathcal{O}_2:
1+J(t,x)
> 0. \]
Then
\begin{enumerate}
\item $\forall\beta \in\mathbf{N}^{\times}:\;\partial^{\beta}\exp h = 
\left(\underset{a \in \mathcal{R}(\beta)}{\sum}\mathcal{P}_{\beta}(a)\;
\underset{j=1}{\overset{\beta}{\prod}}(\partial^j h)^{a_j} \right) 
\exp h .$
\item $\forall \beta \in \mathbf{N}^\times\;\;\;\forall i\in \{1,2\}:$
\[ \partial_i^{\beta}\log (1+J)=- \underset{a\in \mathcal{R}
(\beta)}{\sum}(-1)^{l_{\beta}(a)}\frac{(l_{\beta}(a)-1)!}{(1+
J)^{l_{\beta}(a)}}\; \mathcal{P}_{\beta}(a)\;\underset{j=1}{\overset
{\beta}{\prod}}(\partial_i^j J)^{a_j} . \]
\item $\forall \beta,\gamma \in \mathbf{N}^{\times}:
\;\partial_1^{\gamma}\partial_2^{\beta}\log (1+J)=$
\[ - \underset{b \in \mathcal{R}(\gamma)}{\sum}\;
\underset{a \in \mathcal
{R}(\beta)}{\sum}(-1)^{l_{\beta}(a)+l_{\gamma}(b)}\;\;
\frac{\left( l_{\beta}
(a)+l_{\gamma}(b)-1\right)!}{(1+J)^{l_{\beta}(a)+l_{\gamma}
(b)}}\mathcal{P}_{\beta}(a)\mathcal{P}_{\gamma}(b) \;\;\times \]
\[  \underset{j=1}{\overset{\beta}{\prod}}\;\;\underset{k=1}{\overset{
\gamma}{\prod}}\;(\partial_1^{j}J)^{a_j}(\partial_2^{k}
J)^{b_k} - \underset{a \in \mathcal{R}(\beta)}{\sum}\;\;
\underset{b \in \mathcal{T}(\beta,\gamma)}{\sum}(-1)^{l_{\beta}(a)}\;\;
\frac{(l_{\beta}(a)-1)!}{(1+J)^{l_{\beta}(a)}}\;\;\times \]
\[ \mathcal{P}_{\beta}(a)\mathcal{N}_{\beta,\gamma}(b)
\underset{c \in \underset{i=1}{\overset{\beta}{X}}\mathcal{R}
(b_i,a_i)}{\sum}
\underset{k=1}{\overset{\beta}{\prod}}\;\mathcal{W}_{b_k,a_k}(c^k)
\underset{j=1}{\overset{b_k}{\prod}}(\partial_1^j \partial_2^k 
J)^{c_j^k}. \]
\end{enumerate}
\end{Lem}
\begin{Rem}
One knows that
\[ \forall N \in \mathbf{N}\;\;\;\forall k \in \mathbf{N}^{\times}:\;\;
\mid\{x \in \mathbf{N}^k : \underset{j=1}{\overset{k}{\sum}}x_j = N\}
\mid = \left( \begin{array}{c} N+k-1\\k-1 \end{array} \right). \]
Thus
\[ \forall\;\beta \in \mathbf{N}^{\times}\;:\;\;\mid \mathbf{R}({\beta})
\mid \;\leq \underset{j=1}{\overset{{\beta}}{\sum}}\left( 
\begin{array}{c} {\beta} + k - 1 \\ k - 1 \end{array} \right) = 
\left( \begin{array}{c}2 \beta \\\beta - 1  \end{array} \right)
 < \infty. \]
Therefore, all sums in Lemma 3 are finite. 
\end{Rem}
We skip a simple but lenghty proof of the Lemma 3.
\begin{Lem}
Let $\;(X,d_X)\;,\;(Z,d_Z)\;$ be some metric spaces and let $(Y,p_Y)$ be
a semimetric space. If $X \supseteq K$ is compact and $Q:K \times Y 
\rightarrow Z$ is a continuous map, then
\[ \forall \; \epsilon > 0 \;\;\;\forall \;y_0 \in Y\;\;\;
\exists\;\delta > 0 :
\]
\[ \left( y \in K(y_0,\delta)\right) \Longrightarrow \left(
\forall\;x \in K :d_Z(Q(x,y),
Q(x,y_0)) < \epsilon\right).  \]
\end{Lem}
Proof of this Lemma results from the proof of Lemma IX.3.1 of ~\cite{3}
.\\
Proof of Theorem 1:\\
\emph{Step 1}. Let $(f_n)_{n=0}^{\infty}$ and $ (g_n)_{n=0}^{\infty}$ 
are two sequences convergent in $C^{\infty}(\mathbf{R})$ to $f$ and $g$
, correspondingly. Then,
\[ \forall (\mu,\nu) \in \mathbf{N}^{\times} \times \mathbf{N}\;\;\;\;
\exists c_{\mu \nu} \in \mathbf{R}\;\;\; \;\forall n \in \mathbf{N}:
p_{\mu \nu}(f_n) \leq c_{\mu \nu}.  \]
Let us fix $\;c_{\mu \nu}\;$ for $\;(\mu ,\nu )\in \mathbf{N}^{\times}
\times \mathbf{N} \;$ \\(e.g., $\;\;
 c_{\mu \nu}:= \inf \left\{\tilde{c}_{\mu \nu}
\in \mathbf{R}:\forall \;n \in \mathbf{N}\;\;p_{\mu \nu}(f_n) \leq 
\tilde{c}_{\mu \nu}\right\}) $ \\
and denote 
\[ M_1:= \max \left\{\left( \begin{array}{c} \beta \\ \gamma \end{array}
\right)c_{\alpha \gamma}\;:\;\gamma \in\{0,1, \ldots ,\beta \} 
\right\},\]
\[ M_2 := \max \left\{ \left( \begin{array}{c} \beta \\ \gamma 
\end{array} \right)p_{\alpha (\gamma - \beta)}(g)\;:\;\gamma \in 
\{0,1, \ldots ,\beta \} \right\},  \] 
\[ M_{\alpha \beta}:= \max \{M_1,M_2 \}.\]
Then,
\begin{eqnarray*}
\lefteqn{p_{\alpha \beta}(f_n g_n - f g)\leq}\\
& & \underset{\gamma =0}{\overset {\beta}{\sum}}\left( \begin{array}{c}
\beta \\ \gamma \end{array} \right) p_{\alpha \gamma}(f_n)p_{\alpha 
(\beta - \gamma)}(g_n - g) + \underset{\gamma =0}{\overset{\beta}{\sum}}
\left( \begin{array}{c} \beta \\ \gamma \end{array} \right) 
p_{\alpha \gamma}
(f_n - f) p_{\alpha (\beta - \gamma)}(g) \leq \\
& & M_{\alpha \beta}\;
\underset{\gamma =0}{\overset{\beta}{\sum}}
\left( p_{\alpha (\beta - \gamma)}(g_n - g) + p_{\alpha \gamma}
(f_n - f)\right) \underset{n \rightarrow \infty}
{\longrightarrow} 0 
\end{eqnarray*}
which means that the mapping
\[ C^{\infty}(\mathbf{R})^2 \ni (f,g)
\rightarrow f g \in C^{\infty} (\mathbf{R}) \]
is continuous.

Define $\;h_n:=f_n-f\;$ for $\;n \in \mathbf{N}.\;$ We have $\;e^f-
e^{f_n}=e^f(1-e^{h_n}),\;$ thus
\[ p_{\alpha \beta}(e^f-e^{f_n})\leq \underset{\gamma =0}{\overset
{\beta }
{\sum}}\left( \begin{array}{c}\beta \\\gamma \end{array}\right)
p_{\alpha \gamma}(e^f)p_{\alpha (\beta -\gamma)}(1-e^{h_n}). \]
We get
\[ p_{\alpha \beta}(e^f-e^{f_n})\leq \tilde{M}\;\underset{\gamma =0}
{\overset{\beta}{\sum}}p_{\alpha (\beta - \gamma )}(1-e^{h_n}),\]
where
\[\tilde{M}:=\max \left \{ \left ( \begin{array}{c}\beta \\ \gamma 
\end{array} \right)p_{\alpha \gamma}(e^f):\gamma \in \{0,1, 
\ldots ,\beta
 \} \right \}\in \mathbf{R}. \]
For $\;\gamma =\beta\;$ we have
\[ p_{\alpha 0}(1-e^{h_n})\leq e^{p_{\alpha 0}(h_n)} - 1 \underset{n 
\rightarrow \infty }{\longrightarrow 0}. \]
For $\;\gamma > \beta \geq 0\;$ denote $\;\beta -\gamma := \rho +1\;\;
( \mbox{so}\;\; \rho \geq 0 ),\;$ then
\begin{equation}
p_{\alpha (\rho +1)}(1-e^{h_n}) \leq \underset{\mu =0}{\overset{\rho}
{\sum}}\left( \begin{array}{c}\rho \\\mu \end{array} \right)p_{\alpha 
(\mu +1)}(h_n)p_{\alpha (\rho -\mu )}(e^{h_n}).
\end{equation}
By Lemma 3  (see \emph{1.}) we have
\[\exists M_0 \in \mathbf{R}\;\;\;\forall \mu \in \overline{0,\rho}\;
\;\;\forall \;n \in \mathbf{N}\;:\;\left( \begin{array}{c}\rho \\ \mu
\end{array} \right)p_{\alpha (\rho - \mu)}(e^{h_n}) < M_0, \]
therefore (39) gives
\[ \forall (\alpha,\beta) \ni \mathbf{N}^{\times}\times \mathbf{N}:
p_{\alpha \beta}(e^f-e^{f_n}) \underset{n \rightarrow \infty }{
\longrightarrow 0}. \]

Since the map
\[ \partial :C^{\infty}(\mathbf{R}) \ni f \longrightarrow \partial f \in
 C^{\infty}(\mathbf{R}) \]
is continuos, the mapping $\;\mathcal{A}\;$ defined by (30) is 
continuous as the composition of continuous maps.

\emph{Step 2.} Suppose $\; u \in C^{\infty}(\mathbf{R})\;$ and $\;
\left( \begin{array}{c}a \\ b \end{array} \right) \in \mathbf{R}^2.\;$ 
There is one and only one function $\;g \in C^{\infty}(\mathbf{R}) \;$
such that
\begin{equation}
\left\{ \begin{array}{lcl}g^{\prime \prime}&=& u g \\g(0)&=&a \\
g^{\prime}(0)&=&b \end{array} \right. .
\end{equation}
Making substitution
\begin{equation}
\Psi : \mathbf{R} \ni s \longrightarrow \Psi (s):= \left( \begin{array}
{c}g(s)\\g^{\prime}(s) \end{array} \right)\in \mathbf{R}^2
\end{equation}
in (40) yields
\begin{equation}
\overset{\cdot}{\Psi} = \left( \begin{array}{cc}0&1\\u&0 \end{array}
\right)\Psi,\;\;\;\;\;\Psi (0)= \left( \begin{array}{c}a\\b \end{array}
\right).
\end{equation}
Let us denote the solution of (42) by $\;\Psi (\;\cdot \;;u),\;$ to 
indicate its dependence on $\;u \in C^{\infty}(\mathbf{R}),\;$ and in 
addition let
\[ A: C^{\infty}(\mathbf{R})\ni u \longrightarrow A(u):= \left(
\begin{array}{cc}0&1\\u&0 \end{array} \right) \in C^{\infty}(\mathbf{R},
M_{2\times 2}(\mathbf{R})), \]
\[ B:C^{\infty}(\mathbf{R})\ni h \longrightarrow B(h):= \left(
\begin{array}{cc}0&0\\h&0 \end{array} \right) \in C^{\infty}(\mathbf{R}
,M_{2 \times 2}(\mathbf{R})). \]
( In the sequal $\;A(u)(s):=A(u(s)),\;B(h)(s):=B(h(s))$.)\\
For $\;(\alpha ,\beta )\in \mathbf{N}^{\times}\times \mathbf{N} \;$ and 
$\;\delta > 0 \;$ we denote
\[ s_{\alpha \beta }(\delta ):= \{ h \in C^{\infty}( \mathbf{R}):\;
\forall \gamma \in \overline{0,\beta}\;\;\;p_{\alpha \gamma}(h) < 
\delta \}. \]
We notice that 
\[ \forall (u,h)\in C^{\infty}(\mathbf{R})\times s_{\alpha 0}(1):
 \underset{t \in K_
{\alpha}}{\sup}\parallel A(u + h)(t)\parallel \leq 
1+ p_{\alpha 0}(u), \]
\[ \forall \;h \in C^{\infty}(\mathbf{R})\;:\;\underset{t \in K_{\alpha}
}{\sup}\parallel B(h)(t)\parallel = p_{\alpha 0 
}(h). \]
Now, let us fix $\;u \in C^{\infty}(\mathbf{R})\;$ and $\;\alpha \in 
\mathbf{N}^{\times},\;$ and let us denote
\[ M_{u}:=  4 + p_{\alpha 0}(u). \]

We choose $\;\tau \in \;]\;0,\frac{1}{1+M_{u}}\;[,\;\;\;N(\tau):= 
\min 
\{n \in\mathbf{N}^{\times} : n\tau \geq \alpha \}\;$ and $\;\; K(\tau )
:= [- \tau N(\tau ), \tau N(\tau )].\;\;$ Let $\;N\in \mathbf{N}^
{\times}\;$ is such that $\;N \leq \frac{1}{\tau} \leq N+1 \;$ and let 
$\;d:= \frac{1}{N+2}.\;$ Taking into account that
\[ \forall \;t \in \mathbf{R} \;:\;\Psi (t;u)=\left( 
\begin{array}{c}a\\b \end{array} \right)+
\int _{0}^{t}\;A(u(s))\Psi(s;u)ds \]
we get
\[ \forall t\in[0,\tau]:\;\parallel \Psi(t;u+h)-\Psi (t;u)
\parallel \leq \]
\begin{equation}
 \mid t\mid M_u \underset{s\in[0,\tau]}{\sup}\parallel \Psi (s;
u+h)-\Psi(s;u)\parallel + \mid t\mid \underset{s \in K_{\alpha}}{\sup}
\parallel B(h(s))\Psi (s;u)\parallel .
\end{equation}
Let us fix $\;\epsilon > 0.\;$ Applying Lemma 4 to the map 
\[ Q_u:\mathbf{R}\times C^{\infty}(\mathbf{R})\ni (s,h) 
\longrightarrow Q_u (s,h):= B(h(s))\Psi (s;u)\in \mathbf{R}, \]
and $\;(Y,p_Y)=(C^{\infty}(\mathbf{R}),p_{\alpha 0})\;$ gives 
\[ \exists\delta_0 > 0\;\;\forall h \in s_{\alpha 0}(\delta_0):
\underset{s\in K(\tau )}{\sup }\parallel B(h(s))\Psi(s;u)\parallel < 
\epsilon d^{N(\tau )}. \]
Making use of this in (43) yields
\begin{eqnarray*}
\lefteqn{\forall \;h\in s_{\alpha 0}(\delta _0):\underset{t\in [0,
\tau]}{\sup}\parallel \Psi (t;u+h)-\Psi(t;u)\parallel \leq }\\
& & \tau M_u \underset{s\in [0,\tau]}{\sup}\parallel \Psi (s;u+h)-
\Psi (s;u)\parallel + \epsilon d^{N(\tau)}.
\end{eqnarray*}
Since $\;\tau\;$ is such that $\;0 < \tau < 1-\tau M_u\;$ we get
\[ \forall \;h\in s_{\alpha 0}(\delta _0): \underset{t\in [0,\tau ]}
{\sup}\;\parallel \Psi (t;u+h)-\Psi (t;u)\parallel \leq \epsilon 
d^{N(\tau)}. \]
Applying previous considerations to the case $\;t\in [\tau,2\tau]\;$ 
and making use of
\[ \forall\;h \in s_{\alpha 0}(\delta _0):\;\parallel\Psi (\tau;u+h) - 
\Psi (\tau;u)\parallel \leq \epsilon d^{N(\tau)} \]
one gets
\[ \forall \;h \in s_{\alpha 0}(\delta _0):\underset{s \in [\tau,2\tau 
]}{\sup }\parallel \Psi(s;u+h)-\Psi (s;u)\parallel \leq \epsilon 
d^{N(\tau)-1}. \]
Repeated application of this procedure to $\;[2\tau ,3\tau], \ldots ,
(N(\tau) -1)\tau ,N(\tau)\tau ]\;$ gives 
\begin{eqnarray*}
\lefteqn{ \forall h \in s_{\alpha 0}(\delta _0)\;\;
\forall N \in \overline{1,N(\tau)}\;:}\\ 
& & \underset{s \in [(N-1)\tau ,N\tau]}{\sup}\parallel\Psi 
(s;u+h)-\Psi (s;u)\parallel \leq \epsilon d^{N(\tau)-(N-1)}.
\end{eqnarray*}
Similar reasoning applied to $\;[-\tau ,0],[-2\tau ,-\tau ], \ldots ,
[- N(\tau)\tau,-(N(\tau)-1)\tau]\;$ enables to write
\[ \forall \;h \in s_{\alpha 0}(\delta _0): \underset{s \in K(\tau)}
{\sup}\parallel \Psi(s;u+h)-\Psi(s;u)\parallel \leq \epsilon .\]
Since $\;K_{\alpha}\subseteq K(\tau),\;$ we get
\begin{equation}
\forall \epsilon >0\;\;\exists \delta_0 > 0 \;\;\forall h 
\in s_{\alpha 0}(\delta _0):r_{\alpha
 0}(\Psi(\cdot\;;u+h)-\Psi(\cdot \;;u)) \leq \epsilon .
 \end{equation}
Making use of Lemma 4, Eq. (44) and the estimate
\begin{eqnarray*}
\lefteqn{\forall \;h \in s_{\alpha 0}(1): \;\parallel(\partial \Psi)(t;
u+h)-(\partial \Psi)(t;u)\parallel \leq }\\
& & M_u \parallel \Psi(t;u+h)-\Psi (t;u)\parallel + \parallel B(h(t))
\Psi (t;u)\parallel 
\end{eqnarray*}
we get
\[ \forall \epsilon >0\;\;\exists \delta_1 >0\;\;\forall h \in 
s_{\alpha 1}(\delta_1):r_{\alpha 1}(\Psi(\cdot\;;u+h)-\Psi(\cdot\;;u))
\leq \epsilon. \]
Now, let us consider the identity
\[ \forall \beta \in \mathbf{N}^{\times}:\partial ^{\beta}\left(
\Psi(\cdot
\;;u+h)-\Psi (\cdot \;;u)\right) = \]
\[ \underset{\rho =0}{\overset{\beta -1}{\sum}}\left( 
\begin{array}{c}\beta -1\\ 
 \rho  \end{array}\right)
\left[ A(\partial ^{\rho}(u+h))\partial ^{
\beta -1-\rho}(\Psi(\cdot \;;u+h)-\Psi(\cdot \;;u))\right]-  \] 
\begin{equation}
\underset{\rho =0}{\overset{\beta -1}{\sum}}\left( \begin{array}{c}
\beta -1 \\ \rho \end{array} \right)   B(\partial ^{\rho}
h)\partial ^{\beta -1-\rho}\Psi(\cdot\;;u).
\end{equation}
We notice that there are derivatives of order $\;0,1, \ldots ,\beta -1
\;$ in the right hand side of (45).  
Using Lemma 4 for the map 
\[ Q_u^{\rho}:\mathbf{R}\times C^{\infty} (\mathbf{R})
\ni (s,h)\longrightarrow Q_u^{\rho}(s,h):= B(
\partial ^{\rho}h)\partial ^{\beta -1-\rho}\Psi (s;u)\in 
\mathbf{R}^2, \]
where $\;\rho =0,1, \ldots ,\beta -1\;$ and where $\;(Y,p_Y)=(C^{\infty
}(\mathbf{R}),p_{\alpha \rho}),$\\
we get (by induction)
\[ \forall \beta \in \mathbf{N}\;\;\forall \epsilon >0\;\;
\exists \delta _{\beta}
>0\;\;\forall h \in s_{\alpha \beta}(\delta _{\beta}):r_{\alpha \beta}
(\Psi(\;\cdot \;;u+h)-\Psi(\;\cdot \;;u))\leq \epsilon. \]
Since our considerations apply to any $\;\alpha \in \mathbf{N}^{\times}
\;$ and any $\;u \in C^{\infty}(\mathbf{R})\;$ we obtain that 
\begin{equation}
 C^{\infty}(\mathbf{R})\ni u\longrightarrow \Psi(\;\cdot\;;u)\in C^
{\infty}(\mathbf{R},\mathbf{R}^2) 
\end{equation}
( where $\;\Psi(\;\cdot \;;u)\;$ is a solution of (42) )\\
is a continuous mapping.
Taking into account that
\[ \forall \Phi \in C^{\infty}(\mathbf{R},\mathbf{R}^2)\;\;\;\forall 
(\alpha ,\beta)\in \mathbf{N}^{\times}\times \mathbf{N}:r_{\alpha \beta}
(\Phi)\geq p_{\alpha \beta}(\Phi _1), \]
( where $\Phi (t)=:\left(\begin{array}{c}\Phi _1(t) \\ \Phi _2(t) 
\end{array}\right)$ )
and solving (46) for $\;\left(\begin{array}{c}a\\b \end{array}\right) = 
\left(\begin{array}{c}1\\0 \end{array}\right)\;$ and $\;\left(
\begin{array}{c}a\\b
\end{array} \right) = \left( \begin{array}{c}0\\1 \end{array} 
\right)\;$ gives 
the conclusion that $\;\mathcal{B},\;$ defined by (31), is a continuous 
mapping.\\
\emph{Step 3.} It is clear that both mappings $\;\mathcal{C}$  
and 
$\mathcal{D}$ are continuous. The continuity of $
\mathcal{E}$  can be proved by analogy to the case of 
$\mathcal{A}$ mapping. The continuity of $\mathcal{G}$ 
mapping results from the continuity of the mapping
\[ \forall c \in \{-1,1\}\;\;\;\omega_c : C^{\infty}(\mathbf{R}) \ni f
\longrightarrow \omega_c(f)\in C^{\infty}(\mathbf{R}^2), \]
where
\[ \omega_c(f):\mathbf{R}^2 \ni (t,x) \longrightarrow \omega_c(f)(t,x)
:=f(x+ct)\in \mathbf{R}. \]
The mapping $\mathcal{H}$  is continuous since it is a 
composition of continuous mappings.\\
\emph{Step 4.} What is left is to prove that the mapping $\mathcal{N}
$ is continuous. Suppose $\;(G_n)_{n=0}^{\infty}\;$ is a 
sequence of elements of $\;C_{+}^{\infty}(\mathbf{R}^2)\;$ convergent 
to $\;G\in C_{+}^{\infty}(\mathbf{R}^2).\;$ Denote $\;H_n:= G_n-G. \;$
The sequence $\;(H_n)_{n=0}^{\infty}\;$ converges to zero in $\;C^{
\infty}(\mathbf{R}^2).\;$ Since $\;\forall n \in \mathbf{N}:1+\frac
{H_n}{G} > 0,\;$ we have
\[ \log G_n -\log G = \log\frac{G_n}{G}=\log(1+\frac{H_n}{G}). \]
As $\;H_n \underset{n\rightarrow \infty}{\longrightarrow}0\;\;\mbox{in}
\;\;C^{\infty}(\mathbf{R}^2)\;$ we have
\[ \forall \alpha \in \mathbf{N}^{\times}\;\;\;\exists N_{\alpha}\in 
\mathbf{R}\;\;\;\forall n>N_{\alpha}: q_{\alpha (0,0)}\left(
\frac{H_n}{G}
\right) < \frac{1}{2}. \]
Since $\;\forall x \in \;]\;-1/2,\infty\;[\;\;:\;\;\mid \log(1+x)
\mid \leq 
2\mid x \mid ,\;$ we conclude that\\
$\forall \alpha \in\mathbf{N}^{\times}\;\;\;\exists N_{\alpha}\in 
\mathbf{N}\;\;\;\forall n>N_{\alpha}:$
\[ q_{\alpha (0,0)}\left(\log (1+\frac{H_n}
{G})\right)\leq 2q_{\alpha (0,0)}(\frac{1}{G})q_{\alpha (0,0)}(H_n).\]
Now, suppose $\;\mid \beta \;\mid = \beta _1 + \beta _2 > 0\;$ and 
denote $\;J_n :=\frac{H_n}{G}\;$ for $\;n \in \mathbf{N}.\;$
\\ Since
\[\partial ^{(\beta _1 \beta _2)}J_n = \underset{\gamma _1 =
0}{\overset{\beta _1}{\sum}}\;\underset{\gamma _2 =0}{\overset{\beta _
2}{\sum}}\left( \begin{array}{c}\beta _1\\\gamma _1 \end{array}\right)
\left( \begin{array}{c}\beta_2 \\ \gamma _2 \end{array} \right)
\left(\partial ^{(\gamma _1\gamma _2)}\frac{1}{G} \right) \left( 
\partial ^{(\beta _1-\gamma _1 ,\beta _2-\gamma _2)} H_n \right)\]
we have
\[ q_{\alpha (\beta _1,\beta _2)}(J_n)\leq M_{\alpha 
(\beta _1,\beta _2)}\underset{\gamma _1 =0}{\overset{\beta _1}{\sum}}\;
\underset{\gamma _2=0}{\overset{\beta _2}{\sum}}q_{\alpha (\gamma _1,
\gamma _2)}(H_n), \]
where
\[ M_{\alpha (\beta _1,\beta _2)}:= \max \left\{\left( \begin{array}{c}
\beta _1\\\gamma _1 \end{array} \right)\left( \begin{array}{c}\beta _2
\\\gamma _2 \end{array} \right)q_{\alpha (\gamma _1,\gamma _2)}(\frac
{1}{G}):0\leq \gamma _1\leq \beta _1,0\leq \gamma _2\leq \beta _2
\right\}.  \]
Thus, we see that $\;J_n\underset{n
\rightarrow \infty}{\longrightarrow 0}\;$ in $\;C^{\infty}(\mathbf{R}
^2).\;$ In particular we have
\begin{equation}
 \forall \alpha \in \mathbf{N}^{\times}\;\;\;\exists N^{\alpha}\in 
\mathbf{N}\;\;\;\forall n>N^{\alpha} : \underset{(t,x)\in K_{\alpha}^2}
{\sup}\;\mid \frac{1}{1+J_n(t,x)}\mid \leq 2. 
\end{equation}
For $\;\beta \in \mathbf{N}^2\;$ let $\;D_{\beta}:= \{(i,j)\in 
\mathbf{N}^2:1\leq i+j\leq \mid \beta \mid \},\;\;d_{\beta}:=\mid D
_{\beta}\mid $. \\
By (47) and Lemma 3 (see \emph{2.} and \emph{3.}) we get that for 
$\;\beta \in \mathbf{N}^2,\;\;\mid \beta \mid \geq 1\;$ and $\;\alpha 
\in \mathbf{N}^{\times}\;$ there exists a polynomial $\;Q_{\alpha 
\beta}\in \mathbf{R}[\;x_1, \ldots ,x_{d_{\beta}}\;]\;$ such that 
$\;Q_{\alpha \beta}(0)=0,\;\;\;deg Q_{\alpha \beta} \leq \mid \beta 
\mid \;$ and
\[ \exists N^{\alpha}\in \mathbf{N}\;\;\;\forall n>N^{\alpha}\;\;\;
q_{\alpha (\beta _1 ,\beta 2)}(\log (1+J_n))\leq  \]
\[ Q_{\alpha \beta}\left( q_{\alpha (1,0)}(J_n), q_{\alpha 
(0,1)}(J_n), \ldots ,q_{\alpha (\mid \beta \mid ,0)}(
{J}_n),q_{\alpha (0,\mid \beta \mid )}(J_n) \right). \]
Since $J_n \rightarrow 0$, it follows that
\[ \forall (\alpha ,\beta )\in \mathbf{N}^{\times}\times \mathbf{N}^2
\;\;\;\forall \epsilon >0\;\;\exists N_{\alpha \beta}\in \mathbf{N}\;\;
\forall n>N_{\alpha \beta} : q_{\alpha \beta}(\log (1+J_n))
< \epsilon. \]
But $\;G \in C_{+}^{\infty}(\mathbf{R}^2)\;$ is an arbitrary function,
 therefore the mapping $\;\mathcal{N}\;$ is continuous. Finally, the 
 mapping $\;\mathcal{I},\;$ defined by (38), is continuous as it is 
 a composition of continuous mappings.\\ 
\emph{Step 5.} Denote $\;\mathcal{F}:= C^{\infty}(\mathbf{R})^2.\;$ 
It is clear that the mappings
\[ \mathcal{N}:C^{\infty}(\mathbf{R}^2)\ni F \longrightarrow 
\mathcal{N}(F):=(F(0,\cdot),(\partial_tF)(0,\cdot))\in \mathcal{F} \]
and 
\[ \mathcal{S}:=\mathcal{N}_{{\mid_{\mathcal{M}}}}\] 
are continuous. \\
Corrolary 3 means that
\[ \mathcal{I}\cdot \mathcal{S} = id_{\mathcal{M}}\;\;\;\;\;and \;\;\;
\;\;\mathcal{S}\cdot \mathcal{I} = id_{\mathcal{F}} \]
This completes the proof. 

\section{Concluding Remarks}

Homeomorphism of the set of smooth solutions, $\mathcal{M}\subset 
C^{\infty}(\mathbf{R}^2,\mathbf{R})$,  and 
the space of smooth initial data, $\mathcal{F}$, gives $\mathcal{M}$ 
the structure of a topological manifold modelled on the Fr\'{e}chet
space $\mathcal{F}$.  

A starting point in the geometric quantization of a mechanical system 
is to express the evolution of the system on the phase space in terms 
of  \emph{symplectic} geometry. Can one follow this method in the case 
of the Liouville field theory? We hope to  answer  this 
question in near future.

\section{Acknowledgements}

We are grateful to George Jorjadze for a fruitful discussion. One of us
(K.B.) is greatly indebted to Anatol Odzijewicz for his kind 
hospitality at the Institute of Physics, Warsaw University Branch in 
Bia{\l}ystok, where a part of this paper was done, and wishes to thank 
the So{\l}tan Institute for Nuclear Studies for financial support.

\end{document}